\documentclass[10pt,reqno]{amsart}
\usepackage{amssymb}
\usepackage{amsthm}

\usepackage[foot]{amsaddr}

\usepackage{amssymb,amsmath,amsthm,amsxtra,amsfonts}
\usepackage[margin=4cm]{geometry}

\usepackage{pdfsync}

\usepackage{tikz}

\usepackage{titlesec}
\titleformat{\section}[hang]{\normalfont\large}{\bfseries\thesection.}{.5em}{\bfseries}
\titlespacing{\section}{0pt}{12pt plus 6pt}{3pt}
\titleformat{\subsection}[runin]{\normalfont}{\bfseries\thesubsection.}{.5em}{\bfseries}[.\quad]
\titleformat{\subsubsection}[runin]{\normalfont}{\bfseries\thesubsection.}{.5em}{\bfseries}[.\quad]
\makeatletter
\renewcommand\l@section{\@dottedtocline{1}{0ex}{2em}}
\makeatother

\usepackage[format=plain,labelfont=it,font=normalsize,skip=6pt]{caption}
\DeclareCaptionLabelSeparator{dot}{.~} 
\DeclareCaptionFormat{LabRCaptC}{\hfill#1#2#3}

\usepackage[unicode=true]{hyperref}
\hypersetup{
pdfencoding=auto
,pdftitle={
},pdfauthor={  Alexander I. Aptekarev and Rostyslav Kozhan
}
}

\newcommand{\bbN}{{\mathbb{N}}}

\newcommand{\bbR}{{\mathbb{R}}}

\newcommand{\bbZ}{{\mathbb{Z}}}

\newcommand{\calJ}{{\mathcal J}}

\newcommand{\supp}{\text{\rm{supp}}}

\newcommand{\beq}{\begin{equation}}
\newcommand{\eeq}{\end{equation}}
\newcommand{\ba}{\begin{align*}}
\newcommand{\ea}{\end{align*}}
\newcommand{\veps}{\varepsilon}

\allowdisplaybreaks
\numberwithin{equation}{section}

\newtheorem{theorem}{Theorem}
\newtheorem{lemma}{Lemma}

\newtheorem*{corollary}{Corollary}
\newtheorem*{definition}{Definition}

\theoremstyle{remark}

\theoremstyle{remark}

\begin{document}
\title[Differential equations for recurrence coefficients limits for MOP]{Differential Equations for the Recurrence Coefficients Limits for Multiple Orthogonal Polynomials from a Nevai Class}
\author{Alexander I. Aptekarev$^1$}
\address[$^1$]{
Keldysh Institute of Applied Mathematics, \\
Russian Academy of Science, \\ Miusskaya Pl.4, Moscow 125047,
Russian Federation; e-mail: aptekaa@keldysh.ru
}
\author{Rostyslav Kozhan$^{2}$}
\address[$^2$]{
         Uppsala University\\
         Box 480, 751 06 Uppsala, Sweden; e-mail: kozhan@math.uu.se}
\date{\today}
\begin{abstract}
A limiting property of the nearest-neighbor recurrence coefficients
for multiple orthogonal polynomials from a Nevai class is investigated.
Namely, assuming that the nearest-neighbor coefficients have a limit along rays
of the lattice, we describe it in terms of the solution of a system of
partial differential equations.

In the case of two orthogonality measures the differential equation becomes ordinary. For Angelesco systems, the result is illustrated numerically.
\end{abstract}
\maketitle

\section{Introduction}\label{Intr}
\subsection{Orthogonal polynomials on the real line and the Jacobi matrices}\label{s:OP}
Given a probability measure $\mu$ on $\bbR$ with infinite support, the sequence of its monic orthogonal polynomials $\{P_k\}_{k=0}^\infty$ satisfies the well-known three-term recurrence relation
\begin{equation}\label{eq:3term}
xP_n(x) = P_{n+1}(x) + b_{n} P_n(x) + a_{n-1} P_{n-1}(x)
\end{equation}
with $P_{-1}=0$, $ P_0 = 1$, where the recurrence coefficients
$\{a_k,b_k\}_{k=0}^\infty$  satisfy $a_k>0$, $b_k\in\bbR$.

The corresponding Jacobi matrix is defined to be
\begin{equation}\label{eq:jacobi}
\calJ=\left(
\begin{array}{cccc}
b_0&\sqrt{a_0}&{0}&\\
\sqrt{a_0} &b_1&\sqrt{a_1}&\ddots\\
{0}&\sqrt{a_1} &b_2&\ddots\\
&\ddots&\ddots&\ddots\end{array}\right).
\end{equation}
Assuming $\{a_k\}_{k=0}^\infty$ and $\{b_k\}_{k=0}^\infty$ are bounded, the spectral measure of $\calJ$ with respect to $(1,0,0,\ldots)^T$ coincides with the orthogonality measure $\mu$. Favard's theorem establishes a one-to-one correspondence between all  $\mu$ with compact infinite support and all such bounded self-adjoint Jacobi matrices $\calJ$.

We say that a probability measure on $\bbR$ belongs to the \textit{Nevai class} 
$N(a,b)$ if its Jacobi coefficients (in \eqref{eq:3term}) satisfy $b_n \to b$ and $a_n \to a$ as $n\to\infty$.

Weyl's theorem on compact perturbations implies that any measure in $N(a,b)$ has $\sigma_{ess}(\mu) = [b-2\sqrt{a},b+2\sqrt{a}]$. For the 
converse direction, we have the Denisov--Rakhmanov theorem stating that if
$\sigma_{ess}(\mu) = [\alpha,\beta]$ and $\frac{d\mu}{dx}>0$ a.e. on $[\alpha,\beta]$ then $\mu \in N(\big(\tfrac{\beta-\alpha}{4}\big)^2,\tfrac{\alpha+\beta}{2})$.

See, e.g.,~\cite{Simon11} for more details from the theory of orthogonal polynomials.

\subsection{Multiple orthogonal polynomials and the nearest neighbor recurrence relations}

Let us now describe  multiple orthogonality situation with respect to the vector-measure
$\vec{\mu}:=\{\mu_i\}_{i=1}^d$ on
$\bbR$.
For the rest of the paper we will use the notation $|\vec{v}| := v_1 + ... + v_d $ for any vector-valued object $\vec{v}=(v_j)_{j=1}^d$.

For any $ \vec{n}=(n_1,\ldots,n_d)\in \mathbb{Z}^d_+$, let
$P_{\vec n}$ be the monic polynomial of smallest degree  which
satisfies
\begin{equation}
\label{ortho}
\int P_{\vec{n}}(x)x^k  d \mu_i = 0, \quad k\in\{0,\ldots,n_i-1\}, \quad i\in\{1,\ldots,d\}.
\end{equation}
The polynomial $P_{\vec{n}}(x)$ is called the type II
\emph{multiple orthogonal polynomial} (MOP). Obviously,
$P_{\vec{n}}$ is uniquely determined and $\deg P_{\vec{n}} \leq
|\vec{ n}|$. When $\deg P_{\vec{n}} =  |\vec{
n}|$ the multi-index $\vec{n}$ is said to be \emph{normal}. If all
multi-indices of the lattice $\mathbb{Z}^d_+$ are normal then  the
system of measures  $\{\mu_i\}_{i=1}^d$ is called \emph{perfect}.
It is known  \cite{Ismail2005, VanAssche11}, that (similarly to the
case with one measure) MOPs for the perfect systems satisfy  the
following nearest neighbor recurrence relations (NNRR)
 \begin{equation}
\label{recurrence}
z P_{\vec n}(z) = P_{\vec n+\vec e_j}(z) + b_{\vec n,j}P_{\vec n}(z) + \sum_{i=1}^d a_{\vec n,i}P_{\vec n-\vec e_i}(z),
\end{equation}
where $\vec{e}_j$ is the $j$-th standard basis vector of $\bbR^d$.
Here we have $d$ recurrence relations for $j=1,\ldots,d$. Thus for each $\vec{n}\in \mathbb{Z}^d_+$ we have two sets of the  coefficients for NNRR, namely $\{b_{\vec n,j}\}_{j=1}^d$  and $\{a_{\vec n,i}\}_{i=1}^d$.
Note that for each fixed $j$,  $\{a_{k\vec e_j,j}\}_{k=1}^\infty$ and $\{b_{k\vec e_j,j}\}_{k=1}^\infty$ are the $\{a_k\}_{k=0}^\infty$ and $\{b_k\}_{k=0}^\infty$ from the usual three-term recurrence~\eqref{eq:3term} for the measure $\mu_j$.

In order to define by means of \eqref{recurrence}  the polynomials $\{P_{\vec n}(z)\}$ in unique way the NNRR coefficients cannot be taken  arbitrary.
As was shown in~\cite{VanAssche11},
the recurrence coefficients must  satisfy the compatibility conditions (CC):
 \begin{align}
 \label{Nabla1A}
& b_{\vec{n}+\vec{e}_j,i}-b_{\vec{n},i}=b_{\vec{n}+\vec{e}_i,j}-b_{\vec{n},j},
\quad i<j,
\\
 \label{Nabla2A}
& \det \begin{pmatrix} b_{\vec{n}+\vec{e}_j,i} & b_{\vec{n},i} \\ b_{\vec{n}+\vec{e}_i,j} & b_{\vec{n},j}\end{pmatrix}
=
\sum_{k=1}^d a_{\vec{n}+\vec{e}_j,k}
-
\sum_{k=1}^d a_{\vec{n}+\vec{e}_i,k},
\quad i<j,
\\
 \label{Nabla3A}
 &
 \frac{a_{\vec{n},i}}{a_{\vec{n}+\vec{e}_j,i}}
 =
 \frac{b_{\vec{n}-\vec{e}_i,j} - b_{\vec{n}-\vec{e}_i,i}}{b_{\vec{n},j}-b_{\vec{n},i}},
 \quad i\ne j.
\end{align}
It is not hard to see that these $2d(d-1)$ equalities can be rewritten as
 \begin{align}
 \label{Nabla1}
& \nabla_{j}b_{\vec{n},i}=\nabla_{i}b_{\vec{n},j},
\quad i<j,
\\
 \label{Nabla2}
& b_{\vec{n},j}\nabla_{i}b_{\vec{n},j} - b_{\vec{n},i}\nabla_{j}b_{\vec{n},i}
=\left\langle(\overrightarrow{\nabla}_{j}-\overrightarrow{\nabla}_{i}),\,\vec{a}_{\vec{n}} \right\rangle,
\quad i<j,
\\
 \label{Nabla3}
 & (\nabla_{i}\ln)a_{\vec{n},j}=(\nabla_{j}\ln)\,(b_{\vec{n}-\vec e_{j},i}-b_{\vec{n}-\vec e_{j},j}),
 \quad i\ne j,
\end{align}
where we denote
$$
\nabla_{j}b_{\vec{n},i}:=b_{\vec{n}+\vec e_{j},i}-b_{\vec{n},i},
\quad \overrightarrow{\nabla}_{i}:= (\nabla_{i}, \cdots, \nabla_{i}),
\quad(\nabla_{i}\ln)a_{\vec{n},j}:=
\left(\frac{a_{\vec{n}+\vec e_{i},j}}{a_{\vec{n},j}}-1\right).
$$


The system of difference equations \eqref{Nabla1}--\eqref{Nabla3}
together with the marginal conditions
\begin{equation}\label{eq:marginal}
a_{\vec n,j} =0 , \quad \mbox{ whenever } n_j=0,
\end{equation}
is also called \emph{Discrete Integrable System} (DIS)  for details
see \cite{ADVA}. The boundary problem for DIS
\eqref{Nabla1}--\eqref{Nabla3} in $\mathbb{Z}^d_+$ means the
following. Given the boundary data: coefficients of the
$d$-collections of the three-terms recurrence relations,
corresponding to usual orthogonal polynomials with respect to each
$\{\mu_i\}_{i=1}^d$ measure. Then solving equations
\eqref{Nabla1}--\eqref{Nabla3} we have to find all NNRR
coefficients $\{b_{\vec n,j}\}_{j=1}^d$  and $\{a_{\vec
n,i}\}_{i=1}^d$.

\subsection{Zero asymptotics  and limits of the recurrence coefficients}\label{susec1.3}
Our goal  is to investigate the asymptotic behavior
of the recurrence coefficients $\big\{a_{\vec n,i},b_{\vec
n,i}\big\}$ as $|\:\vec n\:|$ grows.
This behavior is intimately connected to the \emph{asymptotic zero distribution} of multiple orthogonal polynomials $P_{\vec{n}}$.
 To state the problem, we need to
place some restrictions on the way $|\:\vec n\:|$ approaches
infinity as well as the measures $\mu_i$. At the same time we have to be in the class of the perfect systems to keep NNRR. 

 The important example of a perfect system of measures $ \{\mu_i\}$  is the so-called
\emph{Angelesco system} defined by \footnote{If supports of
measures are intervals with nonintersecting interiors then system
$\{\mu_i\}$ is perfect as well.}
\begin{equation}
\label{angelesco}
\supp(\mu_i) = [\alpha_i,\beta_i] ,\quad 
\text{with } \quad \alpha_i<\beta_i<\alpha_{i+1} \quad \text{for all } i.
\end{equation}
 Multiple orthogonal polynomial with respect to  Angelesco system  has the form:
$$
P_{\vec n}(z) =:\displaystyle\prod_{i=1}^{d}\prod_{l=1}^{n_i}(z-x_{\vec{n},i,l}),
\quad x_{\vec{n},i,l}\in [\alpha_i,\beta_i].
$$

\noindent
Moreover, we restrict our attention to sequences of multi-indices
such that
\begin{equation}
\label{multi-indices}
n_i = t_i|\:\vec n\:|+o\left(|\:\vec n\:|\right),
\qquad |\:\vec t\:|=1
\end{equation}
for some $\vec{t} \in (0,1)^d$. 
We denote $\lim_\mathcal{N}$ to be the limit as $|\vec n| \to \infty$ along the sequence of multi-indices satisfying~\eqref{multi-indices}.
Asymptotic zero distribution for  $P_{\vec n}(z)$ (or \emph{limiting zero counting measure}):
\begin{equation}
\label{zero-counting}
\omega (x):= \lim_{\mathcal{N}} \,\,\frac{1}{|\:\vec n\:|} \sum_{i=1}^{d}\sum_{l=1}^{n_i}\delta(x-x_{\vec{n},i,l}),
\end{equation}
for Angelesco systems \eqref{angelesco} with $\mu'_i>0$ a.e. on $[\alpha_i,\beta_i]$ in the regime \eqref{multi-indices} was obtained by
Gonchar and Rakhmanov \cite{GR81}.
To state their result we fix  $\vec t$ as in
\eqref{multi-indices}, and denote
\[
M_{\vec t}\big(\{\alpha_i,\beta_i\}_1^d\big):=\big\{\vec\nu=(\nu_1,\ldots,\nu_d):~\nu_i\in M_{t_i}(\alpha_i,\beta_i), ~i\in\{1,\ldots,d\}\big\},
\]
where $M_t(\alpha,\beta)$ is the set of  positive Borel measures
of mass $t$ supported on $[\alpha,\beta]$.
\begin{theorem}[\cite{GR81}]\label{thm:GR81} \textbf{1)}
There
exists the unique vector of measures $\vec\omega\in M_{\vec
t}\big(\{\alpha_i,\beta_i\}_1^d\big)\,:$
\begin{equation}
\label{ExtrPr}
I[\:\vec\omega\:] = \min_{\nu\in M_{\vec t}(\{\alpha_i,\beta_i\}_1^d)} I[\:\vec\nu\:],  \qquad I[\:\vec\nu\:] := \sum_{i=1}^d\bigg(2I[\nu_i] + \sum_{k\neq i}I[\nu_i,\nu_k]\bigg),\end{equation}
where $I[\nu_i]:=I[\nu_i,\nu_i]$ and
$I[\nu_i,\nu_k]:=-\int\int\log|z-x|{d}\nu_i(x){d}\nu_k(z)$.

\medskip

\noindent
\textbf{2)} Moreover, for the limiting counting measure \eqref{zero-counting} it holds: $\omega=|\vec\omega|$.
\end{theorem}

\noindent
An important feature of the case $d>1$ (in comparison with the classic $d=1$) is the fact that measures $\omega_i$ might no longer be supported on the whole
intervals $[\alpha_i,\beta_i]$ (the so-called \emph{pushing effect}), but in
general it holds that
\begin{equation}
\label{supports}
\supp(\omega_i) =  [\alpha_{\vec t,i},\beta_{\vec t,i}] \subseteq [\alpha_i,\beta_i], \qquad i\in\{1,\ldots,d\}.
\end{equation}
Namely the supports of the extremal measures (not the supports of
the multiple orthogonality measures \footnote{For $d=1$ both of these
notions  coincide.} ) define the recurrence coefficients limits.

\bigskip
To describe the asymptotics of the recurrence coefficients, we
shall need a $(d+1)$-sheeted compact Riemann surface, say $\mathfrak{R}_{\vec{t}}$,
that we realize in the following way. Take $d+1$ copies of
$\overline{\mathbb{C}}$. Cut one of them along the union
$\bigcup_{i=1}^d\big[\alpha_{\vec t,i},\beta_{\vec t,i}\big]$, which
henceforth is denoted by $\mathfrak{R}_{\vec{t}}^{(0)}$. Each of the remaining copies
are cut along only one interval $\big[\alpha_{\vec t,i},\beta_{\vec t,i}\big]$
so that no two copies have the same cut and we denote them by
$\mathfrak{R}_{\vec{t}}^{(i)}$. To form $\mathfrak{R}_{\vec{t}}$, take $\mathfrak{R}_{\vec{t}}^{(i)}$ and glue the banks of
the cut $\big[\alpha_{\vec t,i},\beta_{\vec t,i}\big]$ crosswise to the
banks of the corresponding cut on $\mathfrak{R}_{\vec{t}}^{(0)}$. It can be easily
verified that thus constructed Riemann surface has genus 0. Denote
by $\pi$ the natural projection from $\mathfrak{R}_{\vec{t}}$ to $\overline{\mathbb{C}}$. We
also shall employ the notations $\textbf{z}$ for a point on $\mathfrak{R}_{\vec{t}}$ and $z^{(i)}$ for a point on $\mathfrak{R}_{\vec{t}}^{(i)}$
with $\pi(\textbf{z})=\pi(z^{(i)})=z$.

Since $\mathfrak{R}_{\vec{t}}$ has genus zero, one can arbitrarily prescribe zero/pole
multisets of rational functions on $\mathfrak{R}_{\vec{t}}$ as long as the multisets
have the same cardinality.  Hence, we define $\Upsilon_i$,
$i\in\{1,\ldots,d\}$, to be the rational function on $\mathfrak{R}_{\vec{t}}$ with a
simple zero at $\infty^{(0)}$, a simple pole at $\infty^{(i)}$, and
otherwise non-vanishing and finite. We normalize it so that
$\Upsilon_i(z^{(i)})/z\to1$ as $z\to\infty$. Then the following
theorem holds.
\begin{theorem}[\cite{ADY2}]
\label{thm:recurrence} Let $\{\mu_i\}_{i=1}^d$ be a system of
measures satisfying \eqref{angelesco} and such that
\begin{equation}
\label{weights}
d\mu_i(x) = \rho_i(x)dx,
\end{equation}
where $\rho_i$ is holomorphic and non-vanishing in some
neighborhood of $[\alpha_i,\beta_i]$. Further, let $\mathcal N_{\vec
t}=\{\:\vec n\:\}$ be a sequence of multi-indices as in
\eqref{multi-indices} for some $\vec t\in(0,1)^d$.
Then the
recurrence coefficients $\big\{a_{\vec n,j},b_{\vec
n,j}\big\}$ given by \eqref{recurrence} and \eqref{ortho} satisfy
\begin{equation}
\label{limit}
\lim_{\mathcal N_{\vec t}}a_{\vec n,i} =A_{\vec t,i} \quad \text{and} \quad \lim_{\mathcal N_{\vec t}}b_{\vec n,i} =B_{\vec t,i}, \quad i\in\{1,\ldots,d\},
\end{equation}
where $A_{\vec t,i}$ and $B_{\vec t,i}$ are constants:  $\textbf{z}^2\Upsilon_i(z^{(0)}) = A_{\vec t,i}(z+B_{\vec
t,i}) + \mathcal O\big(z^{-1}\big)$ as $z\to\infty$.
\end{theorem}

\noindent
\textbf{Remarks.} \textbf{1)} We note that Theorem~\ref{thm:recurrence} is valid for $d=1$ as well.

\noindent
\textbf{2)} It is not too difficult to extend the proof  (from \cite{GR81}) of Theorem~\ref{thm:GR81} to include the case of touching intervals.

\noindent
\textbf{3)} We also can affirm (at least for $d=2$) that Theorem~\ref{thm:recurrence} remains valid for the case of touching intervals (technicalities can be taken from \cite{AVAY}) and for weight functions \eqref{weights} with singularities  of the types: Jacobi and Fisher-Hartwig weights \cite{Yattselev16}.  $\blacksquare$

\bigskip
%
%

Let us make the following definition by analogy with the scalar case (see Section~\ref{s:OP}).
\begin{definition}
We say that a perfect system of measures $\{\mu_i\}_{i=1}^d$
belongs to the \textit{multiple Nevai class} if for each $ i\in\{1,\ldots,d\}$ the  limits
$$
\lim_{\mathcal N_{\vec t}}a_{\vec n,i} \quad \text{and} \quad \lim_{\mathcal N_{\vec t}}b_{\vec n,i}
$$
exist along each sequence~\eqref{multi-indices} for any $\vec{t}\in [0,1]^d$, $|\vec{t}|=1$.
\end{definition}

Perfect systems from multiple Nevai class appear naturally in various contexts~\cite{A, AKLR,
AKS, K95, VanAssche14}, e.g., in random matrix theory~\cite{DFK}. Note that if a system of measures belongs to a multiple Nevai
class, then the recurrence along the step-line has asymptotically
periodic recurrence coefficients. 


Notice that Theorem~\ref{thm:recurrence} can be viewed as a partial analogue of the Denisov--Rakhmanov theorem,
and Angelesco systems from Theorem~\ref{thm:recurrence}
belong to the multiple Nevai class.
It is an interesting open  problem to generalize this analogue of Denisov--Rakhmanov result 
to more general measures (i.e. to Angelesco systems with $\mu_j>0$ a.e. on $\sigma_{ess}(\mu_j)$).

\medskip

The organization of the paper is as follows. In Section~\ref{s:PDE}
we state and prove our main result: a conditional theorem on
partial differential equations for the limiting value  (in the
regime \eqref{multi-indices}) of the NNRR coefficients. In
Section~\ref{DE_Ang_d2} we discuss the special case of two $d=2$
orthogonality measures when our partial differential equations
become ordinary differential equations. In
Section~\ref{AngParametriz}, using a parametrization of
$\mathfrak{R}_{\vec t}$ from~\cite{LyTu}, we give a constructive
procedure for determination of  limits in~\eqref{limit}. Finally, in Section~\ref{Ang_d2Num}
we present numeric illustrations.

\section{Differential equations for the limits of NNRR coefficients}\label{s:PDE}

\subsection{Construction of the approximating functions}\label{ss:approximation}
For the rest of the paper, let us denote
\begin{equation}\label{eq:S}
S_{d-1} : = \{\vec{s}\in [0,1]^{d-1}: |\vec{s}| \le 1\}.
\end{equation}

Assume that $\{\mu_j\}_{j=1}^d$ form a perfect system from the multiple Nevai class.

This means that there exist $S_{d-1}\to \bbR$ functions $A_j(\vec{s}), B_j(\vec{s})$ ($1\le j \le d$) defined via
\begin{align}
\label{Afunc} & A_j(\vec s) = \lim_{\mathcal{N}} a_{\vec n,j}, \\
\label{Bfunc}  & B_j(\vec s) = \lim_{\mathcal{N}} b_{\vec n,j},
\end{align}
where $\lim_\mathcal{N}$ notation is defined in Section~\ref{susec1.3} with $\vec{t} = \{ \vec{s},1-|\vec s|\}$ (that is, $\vec{s}$ consists of the first $d-1$ coordinates of $\vec{t}$ which defines the direction of the approach to infinity).

In this paper we investigate the possibility of describing functions $\{A_j,B_j\}_{j=1}^d$ through differential equations. This is done in Theorem~\ref{thm:mainMULTI} below.

Before stating the main result, let us introduce the families of approximations $A_j^{(m)}$ and $B_j^{(m)}$ of the limiting functions $A_j$ and $B_j$.

Fix $m\in \bbZ_+$ and $1\le j \le d$. We take all the coefficients $\{a_{\vec n, j}\}$ with $|\vec{n}| = m$ and form an approximating function $A^{(m)}(\vec s)$ as follows. First, for any $\vec n$ with $|\vec{n}| = m$, define $\vec{s}\in S_{d-1}$ via $s_j:= \tfrac{n_j}{m}$ ($1\le j \le d-1$) and let
$$
A^{(m)}_j(\vec{s})= a_{\vec n, j}.
$$
For points in $\tfrac{1}{m} \bbZ_+^{d-1}$ that are not in $S_{d-1}$ we can choose $A^{(m)}$ to be zero.
Then we can extend $A^{(m)}_j$ to the rest of the simplex $S_{d-1}$ via the multilinear interpolation which can be written as follows. Choose a cube $K$ of side length $\tfrac1m$ with vertices in $\tfrac{1}{m} \bbZ_+^{d-1}$; let us denote them $\{P^{(k)},Q^{(k)}\}_{k=1}^{2^{d-2}}$, where for each $k$, vertices $P^{(k)}$ and $Q^{(k)}$ are opposite of each other. If $P^{(k)} = (p^{(k)}_1,\ldots,p^{(k)}_{d-1})$ and $Q^{(k)} = (q^{(k)}_1,\ldots,q^{(k)}_{d-1})$ then we let
\begin{equation}\label{eq:multilinear}
A^{(m)}_j(\vec u):= \sum_{k=1}^{2^{d-2}} \left[ A^{(m)}_j(P^{(k)}) \prod_{l=1}^{d-1} \frac{q^{(k)}_l-u_l}{q_l^{(k)}-p_l^{(k)}}
+
A^{(m)}_j(Q^{(k)}) \prod_{l=1}^{d-1} \frac{p^{(k)}_l-u_l}{p_l^{(k)}-q_l^{(k)}}
\right].
\end{equation}
for $\vec u \in K$.


The main features of this multilinear interpolation function~\eqref{eq:multilinear} that are important to us are:

\medskip

\textbf{1.} The right-hand side of~\eqref{eq:multilinear} agrees with the left-hand side of~\eqref{eq:multilinear} when $\vec{u}\in \{P^{(k)},Q^{(k)}\}_{k=1}^{2^{d-2}}$, so that the function is well defined at the vertices of our cubes;

\smallskip

\textbf{2.} For $\vec{u}$ belonging to any face of a cube $K$, the expression~\eqref{eq:multilinear} reduces to the multilinear interpolation of one dimension lower over the vertices of that face. As a result, ~\eqref{eq:multilinear} on a face of a cube $K$ will agree with~\eqref{eq:multilinear} defined through another cube sharing the same face. So the function $A^{(m)}$ is well-defined on $S_{d-1}$. Moreover, it is continuous on $S_{d-1}$ and is differentiable on the interiors of each of the cubes $K$;

\smallskip

\textbf{3.} In each of the $d-1$ variables $u_l$, the function $A_j^{(m)}$ is linear within each of the cubes $K$. This will be used in the proof of Theorem~\ref{thm:main} below;

\smallskip

\textbf{4.} Partial derivatives of the right-hand side of~\eqref{eq:multilinear} are linear functions along each path parallel to the coordinate axes. In particular, it implies that the maxima and minima over $K$ of partial derivatives of $A^{(m)}_j$ are attained at $\{P^{(k)},Q^{(k)}\}_{k=1}^{2^{d-2}}$. This will be used in the proof of Lemma~\ref{lm:derivativesEst} below.

\smallskip

We can do the same construction with coefficients $b_{\vec n, j}$ to form the multilinear approximations $B^{(m)}_j:S_{d-1}\to \bbR$ for functions $B_j$.

Notice that \eqref{Afunc}--\eqref{Bfunc} implies pointwise convergence  $A^{(m)}_j$ and $B^{(m)}_j$ on $S_{d-1}$ to $A_j$ and $B_j$, respectively, as $m\to\infty$.

\subsection{The main theorem}

For the rest of the paper we assume that the functions $A_j$ and $B_j$ ($1\le j \le d$) are piecewise continuously differentiable on $S_{d-1}$ in the following sense. We suppose that $S_{d-1}$ can be decomposed into a finite union of closed sets $\{D_i\}$ such that:

\medskip

\textbf{(i)}  $A_j$ and $B_j$ are differentiable on the interior $\operatorname{Int}(D_i)$;

\smallskip

\textbf{(ii)} Each of the partial derivatives of $A_j$ and $B_j$ are continuous $\operatorname{Int}(D_i)$ and can be continuously extended to $D_i$.

\smallskip

Note that the latter condition means that each of the partial derivatives of $A_j$ and $B_j$ is uniformly continuous on $\operatorname{Int}(D_i)$, a fact that we use in the proof of Lemma~\ref{lm:derivativesEst}.

We also assume that sets $D_i$ are not pathological, in particular, the closure of $\operatorname{Int}(D_i)$ is assumed to be $D_i$.

Recall that $\{\vec e_j\}_{j=1}^d$ is the standard basis of $\bbR^d$. For the notational convenience, let us denote $\vec{\delta}_j$ ($1\le j \le d-1$) to be the $j$-th standard basis vector in $\bbR^{d-1}$, while $\vec{\delta}_d$ to be the zero vector in $\bbR^{d-1}$.


\begin{theorem}\label{thm:mainMULTI}
Assume that we have a perfect system $\{\mu_j\}_{j=1}^d$ from the multiple Nevai class satisfying the conditions
\begin{itemize}
 \item[(i)] $A_j$ and $B_j$ are piecewise continuously differentiable  
 on $S_{d-1}$ for each $1\le j \le d$;
 \item[(ii)] For each $1\le j \le d$, we have uniform convergence:
\begin{align}\label{fastA}
|A^{(m)}_j(\vec{s}) - A_j(\vec{s}) | & \le  o(\tfrac1m), \\
\label{fastB}
|B^{(m)}_j(\vec{s}) - B_j(\vec{s}) | & \le o(\tfrac1m),
\end{align}
as $m\to\infty$, where sequences $o(\tfrac1m)$ are independent of $\vec s\in S_{d-1}$.
\end{itemize}
Then the limiting functions $A_j$ and $B_j$, $1\le j \le d$, satisfy the following system of $2d(d-1)$ differential equations:
\begin{align}
\label{PDE1}
& \nabla B_{i}(\vec s) \cdot \left(\vec{\delta}_j-\vec{s}\right)
=
\nabla B_{j}(\vec s) \cdot \left(\vec{\delta}_i-\vec{s}\right),\quad i<j,
 \\
\label{PDE2}
& B_j(\vec s)\, \nabla B_{j}(\vec s) \cdot \left(\vec{\delta}_i-\vec{s}\right) -B_i(\vec s) \nabla B_{i}(\vec s) \cdot \left(\vec{\delta}_j-\vec{s}\right)
=
\left( \sum_{l=1}^d \nabla A_l\right) \cdot \left(\vec{\delta}_j-\vec{\delta}_i\right), \quad i<j,
 \\
\label{PDE3}
& A_j(\vec s)\, \nabla \left( B_i(\vec s) -  B_j(\vec s)\right)\cdot \left(\vec{s} - \vec{\delta}_j\right) + (B_i(\vec s)-B_j(\vec s)) \, \nabla A_j(\vec s) \cdot \left(\vec{\delta}_i - \vec{s}\right) = 0, \quad i\ne j.
\end{align}
\end{theorem}
\noindent  In the system~\eqref{PDE1}--\eqref{PDE3}, $\vec{u}\cdot\vec{v}$ stands
for the standard inner product in $\bbR^{d-1}$, and $\nabla$ for
the gradient operator for
a function of $d-1$
variables.

\medskip

\noindent
\textbf{Remarks.} 

\noindent \textbf{1)}
 Condition (i) is fulfilled for Angelesco systems from Theorem~\ref{thm:recurrence}.
 This  follows from
smoothness of the dependence of the residues of $\Upsilon$
on $\vec{t}$. {We show it explicitly for $d = 2$ in the last
section.} As for (ii),
\eqref{fastA}--\eqref{fastB} holds uniformly on compacts of $\operatorname{Int}(S_{d-1})$  (this follows from  the proof of
Theorem~\ref{thm:recurrence}). Whether this can be extended to the whole $S_{d-1}$ is still unknown.

\noindent \textbf{2)} Since the system $\{\mu_j\}_{j=1}^d$ is in
the multiple Nevai class determined by the functions
$\{A_j,B_j\}_{j=1}^d$, each of the measures $\mu_j$ is in the Nevai
class, in particular its essential support is an interval. These
intervals
(together with~\eqref{eq:marginal}) allow one to establish boundary
conditions for the functions $\{A_j,B_j\}_{j=1}^d$. We do this
explicitly for $d=2$ in the next section.

%


\subsection{Convergence of the derivatives}

In order to prove Theorem~\ref{thm:mainMULTI}, we will need to control the derivatives of our approximation functions. This is the purpose of the following lemma.

\begin{lemma}\label{lm:derivativesEst}
Suppose (i)--(ii) of Theorem~\ref{thm:mainMULTI} hold. Then for $1\le k\le d-1$ and any point $\vec{s}_0$ in $\operatorname{Int}(D_i)$, there exists a neighbourhood $U(\vec{s}_0)\subset\operatorname{Int}(D_i)$ such that
\begin{align}
\label{fastAprime}
\left|\frac{\partial}{\partial s_k} A^{(m)}_j(\vec{s}) - \frac{\partial}{\partial s_k} A_j(\vec{s}) \right| & \le  o(1),  \\
\label{fastBprime}
\left|\frac{\partial}{\partial s_k} B^{(m)}_j(\vec{s}) - \frac{\partial}{\partial s_k} B_j(\vec{s}) \right| & \le  o(1),
\end{align}
for all $\vec{s}\in U(\vec{s}_0)$ as $m\to\infty$, where $o(1)$ is independent of $\vec s\in U(\vec{s}_0)$.
\end{lemma}
\noindent
\textbf{Remark.}

\noindent
Partial derivatives of $A^{(m)}_j$ and $B^{(m)}_j$ have jump discontinuities along each side of the $\tfrac1m\bbZ_+^{d-1}$ cubes (see Section~\ref{ss:approximation}). At a point of discontinuity, we interpret $\frac{\partial}{\partial s_k} A^{(m)}_j(\vec{s})$ and $\frac{\partial}{\partial s_k} B^{(m)}_j(\vec{s})$ in~\eqref{fastAprime} and~\eqref{fastBprime} as one of the limiting values of these functions from the inside of one of the cubes.

\begin{proof}
Fix $j$. Let us prove~\eqref{fastAprime} for $k=1$.

 Choose $M_1\in\bbN$ large enough so that a cube with side length $\tfrac2{M_1}$ centered at $\vec{s}_0$ belongs to $D_i$. Let $U(\vec{s}_0)$ be the cube centered at $\vec{s}_0$ of side length $\tfrac{1}{M_1}$.

Let $\veps>0$ be arbitrary. By the discussion in the beginning of the section, $\frac{\partial}{\partial s_1} A_j$ is uniformly continuous on $D_i$. We can therefore find $M_2\in\bbN$ so that
\begin{equation}\label{eq:temp1}
\left| \frac{\partial}{\partial s_1} A_j(\vec{s})  - \frac{\partial}{\partial s_1} A_j(\vec{u})  \right| \le \frac\veps4
\end{equation}
for all $\vec{s}$ and $\vec{u}$ in $D_i$ satisfying $||\vec{s}-\vec{u}|| \le \tfrac1{M_2}$. By ~\eqref{fastA}  we can find $M_3\in\bbN$ so that
\begin{equation}\label{eq:temp2}
m |A^{(m)}_j(\vec{s}) - A_j(\vec{s}) | \le  \frac\veps4
\end{equation}
for all $\vec s\in S_{d-1}$ and $m\ge M_3$. Now let $M= \max\{M_1,M_2,M_3\}$. 

For any $\vec{s}$ in $U(\vec{s}_0)$ and any $m\ge M$, choose a cube $K(m)$ of side length $\tfrac{1}{m}$ containing $\vec{s}$ whose vertices are at $\tfrac1m \bbZ_+^{d-1}$ (as in Section~\ref{ss:approximation}). By the construction, $K$ belongs to $D_i$, and ~\eqref{eq:temp1} and~\eqref{eq:temp2} hold for our $m$.

Let us first show that the inequality~\eqref{fastAprime} holds for the case when $\vec{s}$ is a vertex of $K$. If $\vec{s}+\tfrac1m \vec{e}_1$ is also a vertex of $K$ (arguments for $\vec{s}-\tfrac1m \vec{e}_1$ are identical),  then by the discussion after~\eqref{eq:multilinear},
\begin{multline*}
\left| \frac{\partial}{\partial s_1} A^{(m)}_j(\vec{s}) - \frac{\partial}{\partial s_1} A_j(\vec{s}) \right|
=
\left| m\left[ A^{(m)}_j(\vec{s}+\tfrac1m \vec{e}_1)- A^{(m)}_j(\vec{s})\right] - \frac{\partial}{\partial s_1} A_j(\vec{s}) \right|
\\
\le m \left| (A^{(m)}_j - A_j)(\vec{s}+\tfrac1m \vec{e}_1) \right|
+ m \left| (A^{(m)}_j - A_j)(\vec{s}) \right|
\\
+ \left| m\left[ A_j(\vec{s}+\tfrac1m \vec{e}_1)- A_j(\vec{s})\right] - \frac{\partial}{\partial s_1} A_j(\vec{s})  \right|
\le \frac\veps4+\frac\veps4 + \left| \frac{\partial}{\partial s_1} A_j(\vec{\theta})  - \frac{\partial}{\partial s_1} A_j(\vec{s})  \right|
\end{multline*}
for some $\vec\theta \in (\vec{s},\vec{s}+\tfrac1m \vec{e}_1)$. Here we used~\eqref{eq:temp2} twice and the Mean Value Theorem. The last expression is $\le 3\veps/4$ by~\eqref{eq:temp1}.

Now if $\vec{s}$ is not a vertex of $K$, then by the discussion after~\eqref{eq:multilinear}, there are vertices $\vec{z}_1$ and $\vec{z}_2$ of $K$ such that
$\frac{\partial}{\partial s_1} A^{(m)}_j(\vec{z}_1)\le \frac{\partial}{\partial s_1} A^{(m)}_j(\vec{s}) \le \frac{\partial}{\partial s_1} A^{(m)}_j(\vec{z}_2)$. By~\eqref{eq:temp1}, $\frac{\partial}{\partial s_1} A_j(\vec{z}_2)-\tfrac\veps4 \le \frac{\partial}{\partial s_1} A_j(\vec{s}) \le \frac{\partial}{\partial s_1} A_j(\vec{z}_1)+\tfrac\veps4$. Combining these two inequalities together with the estimate at the vertices, we get $\left| \frac{\partial}{\partial s_1} A^{(m)}_j(\vec{s}) - \frac{\partial}{\partial s_1} A_j(\vec{s}) \right| \le \veps$.
\end{proof}

\subsection{Proof of Theorem~\ref{thm:mainMULTI}}


Let $\vec{s}\in S_{d-1}$ belongs to the interior of some $D_i$. Choose a neighbourhood $U\subset D_i$ of $\vec{s}$ as in Lemma~\ref{lm:derivativesEst}. We can assume $\overline{U}\subset D_i$ (just shrink $U$ if needed). Let a sequence of multi-indices $\vec{n}$ be given satisfying~\eqref{multi-indices} with $\vec{t} = \{\vec{s},1-|\vec{s}|\}$, and as a result~\eqref{Afunc},~\eqref{Bfunc} also. For each such $\vec{n}$, let $m=|\vec{n}|$ and define $\vec s\,^{(m)} \in S_{d-1}$ with $s^{(m)}_j=\frac{n_j}{m}$. 
Then $\vec{s}\,^{(m)} \to \vec{s}$. For each $m$ let $K_m$ be a cube of side length $\tfrac{1}{m+1}$ containing $\vec{s}\,^{(m)}$ whose vertices are at $\tfrac1{m+1} \bbZ_+^{d-1}$ (as in Section~\ref{ss:approximation}). We consider $m$ large enough so that each $K_m$ belongs to $U$.



Let $1\le i \le d-1$. Notice that by Taylor's theorem
\begin{align}
\label{eq:new1}
a_{\vec n+\vec{e}_i,j} & =A_{j}^{(m+1)}(\tfrac{m}{m+1} \vec s\,^{(m)} + \tfrac{1}{m+1}\vec{\delta}_i)
\\
\label{eq:new2}
 & = A_{j}^{(m+1)}(\vec s\,^{(m)}) + \nabla A_{j}^{(m+1)}(\vec s\,^{(m)}) \cdot \left(\tfrac{m}{m+1} \vec s\,^{(m)} + \tfrac{1}{m+1}\vec{\delta}_i-\vec s\,^{(m)}\right) + o(\tfrac{1}{m})
\\
\label{eq:new3}
&= A_{j}^{(m+1)}(\vec s\,^{(m)}) + \tfrac{1}{m+1} \nabla A_{j}^{(m+1)}(\vec s\,^{(m)}) \cdot \left(\vec{\delta}_i-\vec s\,^{(m)}\right) + o(\tfrac{1}{m})
\\
\label{eq:new4}
& = A_{j}(\vec s\,^{(m)}) + \tfrac{1}{m+1} \nabla A_{j}(\vec s\,^{(m)}) \cdot \left(\vec{\delta}_i-\vec s\,^{(m)}\right) + o(\tfrac{1}{m}),
\end{align}
where on the last step we used (ii) of Theorem~\ref{thm:main} and Lemma~\ref{lm:derivativesEst}.
However the $o(\tfrac{1}{m})$ error term in~\eqref{eq:new2} is dependent on $\vec s\,^{(m)}$ and can in principle be non-uniform. To justify uniformity in~\eqref{eq:new4} we proceed as follows.
We start with~\eqref{eq:new1},
and note that $\tfrac{m}{m+1} \vec{s}\,^{(m)} + \tfrac{1}{m+1}\vec{\delta}_i = \vec{s}\,^{(m)}  + \sum_{l=1}^{d-1} \vec\Delta^{(m,i)}_l$ where $\vec\Delta^{(m,i)}_i = \frac{1-s_i^{(m)}}{m+1} \vec\delta_i$ and $\vec\Delta^{(m,i)}_l = - \frac{s_l^{(m)}}{m+1} \vec\delta_l$ for $l\ne i$. These $\vec\Delta^{(m,i)}_l$ are just the increment $\tfrac{1}{m+1} \left(\vec{\delta}_i-\vec s\,^{(m)}\right)$ from $\vec{s}\,^{(m)}$ to $\tfrac{m}{m+1} \vec{s}\,^{(m)} + \tfrac{1}{m+1}\vec{\delta}_i$ separated in coordinates, and $\sum_{l=1}^{d-1} \vec\Delta^{(m,i)}_l = \tfrac{1}{m+1} \left(\vec{\delta}_i-\vec s\,^{(m)}\right)$. Now recall that the multilinear approximation function $A_j^{(m)}$~\eqref{eq:multilinear} is linear along coordinate axes, so applying this for each of the $d-1$ increment we get:
\begin{align*}
a_{\vec n+\vec{e}_i,j} & =A_{j}^{(m+1)}\left(\vec{s}\,^{(m)}  + \sum_{l=1}^{d-1} \vec\Delta^{(m,i)}_l\right)\\
 & = A_{j}^{(m+1)}\left(\vec{s}\,^{(m)}  + \sum_{l=1}^{d-2} \vec\Delta^{(m,i)}_l\right) +
 \nabla A_{j}^{(m+1)}\left(\vec{s}\,^{(m)}  + \sum_{l=1}^{d-2} \vec\Delta^{(m,i)}_l\right) \cdot \vec\Delta^{(m,i)}_{d-1} = \ldots
\\
&= A_{j}^{(m+1)}\left(\vec{s}\,^{(m)} \right) +
 \sum_{p=1}^{d-1} \nabla A_{j}^{(m+1)}\left(\vec{s}\,^{(m)}  + \sum_{l=1}^{p-1} \vec\Delta^{(m,i)}_l\right) \cdot \vec\Delta^{(m,i)}_{p}
\\
&
= A_{j}\left(\vec{s}\,^{(m)} \right) +
 \sum_{p=1}^{d-1} \nabla A_{j}\left(\vec{s}\,^{(m)}  + \sum_{l=1}^{p-1} \vec\Delta^{(m,i)}_l\right) \cdot \vec\Delta^{(m,i)}_{p} + o(\tfrac1m),
\end{align*}
where on the last step we used (ii) of Theorem~\ref{thm:main} and Lemma~\ref{lm:derivativesEst} (notice that now $o(\tfrac1m)$ is uniform!). Now for any $p$, $\nabla A_{j}\left(\vec{s}\,^{(m)}  + \sum_{l=1}^{p-1} \vec\Delta^{(m,i)}_l\right) = \nabla A_{j}\left(\vec{s}\,^{(m)}\right) + o(1)$ (with uniform $o(1)$), since $\vec\Delta^{(m,i)}_l=o(1)$ for each $l$ and $\nabla A_{j}$ is continuous and therefore uniformly continuous on $\overline{U}$. Plugging this into the last equation and using $\vec\Delta^{(m,i)}_p=o(\tfrac1m)$ implies~\eqref{eq:new4} with uniform $o(\tfrac1m)$.

Similar arguments give us for $1\le i \le d-1$,
\begin{equation*}
a_{\vec n-\vec{e}_i,j}  = A_{j}(\vec{s}\,^{(m)}) + \tfrac{1}{m-1} \nabla A_{j}(\vec s\,^{(m)}) \cdot \left(\vec{s}\,^{(m)} - \vec{\delta}_i\right) + o(\tfrac{1}{m})
\end{equation*}
with uniform $o(\tfrac1m)$. For $i=d$, we get the following expressions instead:
\begin{align*}
a_{\vec n+\vec{e}_d,j} & = A_{j}(\vec{s}\,^{(m)}) -\tfrac{1}{m+1} \nabla A_{j}(\vec s\,^{(m)}) \cdot \vec{s}\,^{(m)} + o(\tfrac{1}{m});
\\
a_{\vec n-\vec{e}_d,j} &= A_{j}(\vec{s}\,^{(m)}) + \tfrac{1}{m-1} \nabla A_{j}(\vec s\,^{(m)}) \cdot \vec{s}\,^{(m)} + o(\tfrac{1}{m})
\end{align*}
with uniform $o(\tfrac1m)$.
Notice that these expressions for $a_{\vec n\pm\vec{e}_d,j}$ agree with the expressions for $a_{\vec n\pm\vec{e}_i,j}$ (with $i\le d-1$) if we adopt our notation $\vec{\delta}_d = \vec{0}\in\bbR^{d-1}$.

Analogous equalities hold for the $b$-coefficients and the corresponding $B_j$ functions.

Let us now plug these equalities into~\eqref{Nabla1}. For any $1\le i<j\le d$, we get:
\begin{multline*}
B_{i}(\vec{s}\,^{(m)}) + \tfrac{1}{m+1} \nabla B_{i}(\vec s\,^{(m)}) \cdot \left(\vec{\delta}_j-\vec{s}\,^{(m)}\right) + o(\tfrac{1}{m})
-B_{i}(\vec{s}\,^{(m)})
\\
=
B_{j}(\vec{s}\,^{(m)}) + \tfrac{1}{m+1} \nabla B_{j}(\vec s\,^{(m)}) \cdot \left(\vec{\delta}_i-\vec{s}\,^{(m)}\right) + o(\tfrac{1}{m})
-B_{j}(\vec{s}\,^{(m)}).
\end{multline*}
Now multiply by $m+1$, and take limit as $m\to\infty$.
Using continuity of $\nabla B_j$'s, we obtain ~\eqref{PDE1}.

Similar computations show that ~\eqref{Nabla2} leads to ~\eqref{PDE2} (for any $1\le i<j\le d$)
and~\eqref{Nabla3} produces ~\eqref{PDE3} (for any $i \ne j$, $1\le i,j \le d$).


\section{$d=2$ case: system of ordinary differential equations}\label{DE_Ang_d2}
\subsection{The main theorem: $d=2$}
In the case $d=2$, we have four functions $A_1,A_2,B_1,B_2$ of one variable $s\in[0,1]$, and the corresponding differential system takes the form stated below.

\begin{theorem}\label{thm:main}

$\textbf{1)}$
Assume that we have a perfect system $\mu_1,\mu_2$ from the multiple Nevai class satisfying the conditions
\begin{itemize}
 \item[(i)] $A_j$ and $B_j$ are piecewise continuously differentiable  
 on $[0,1]$ for each $1\le j \le 2$;
 \item[(ii)] For each $1\le j \le 2$, we have uniform convergence:
\begin{align}\label{fastAs}
|A^{(m)}_j(s) - A_j(s) | & \le  o(\tfrac1m), \\
\label{fastBs}
|B^{(m)}_j(s) - B_j(s) | & \le o(\tfrac1m),
\end{align}
as $m\to\infty$, where sequences $o(\tfrac1m)$ are independent of $s\in [0,1]$.
\end{itemize}
Then the limiting functions $A_j$ and $B_j$, $1\le j \le 2$, satisfy the following system of ordinary differential equations:
\begin{equation}\label{3ODE}
\begin{pmatrix}
s B(s) & 0 & (1-s)A_1(s) \\
0 & (1-s) B(s) & s A_2(s) \\
1 & 1 & s(1-s) B(s)
\end{pmatrix}
\begin{pmatrix}
A_1'(s) \\
A_2'(s) \\
B'(s)
\end{pmatrix}
=
\begin{pmatrix}
0 \\
0 \\
0
\end{pmatrix}
\end{equation}
where
\begin{equation}\label{eq:Bs}
B(s) = B_2(s) - B_1(s), \quad B_1'(s) = \frac{A_1'(s)+A_2'(s)}{sB(s)}, \quad B_2'(s) = - \frac{A_1'(s)+A_2'(s)}{(1-s)B(s)}.
\end{equation}

\bigskip

$\textbf{2)}$
Suppose an Angelesco system  satisfies conditions of Theorem~\ref{thm:recurrence} and (ii). Then there exist $c_1, c_2\in (0,1)$
such that the functions $A_1,A_2,B_1,B_2:[0,1]\to\bbR$ are smooth on $[0,c_1)$ and $(c_2,1]$, and satisfy the system of differential equations
\begin{equation}\label{ODEsys}
\begin{cases}
 & (1+s) s C_1'(s) + 4s  C_1(s) + (2-s)(1-s) C_2'(s)   -4(1-s) C_2(s)  =0
\\
 &  \displaystyle\frac{s^2 C_1'(s)}{C_1(s)} =  \frac{(1-s)^2 C_2'(s)}{C_2(s)} -2
\end{cases}
\end{equation}
with initial/boundary conditions
\begin{align}\label{ODEbc}
& \begin{cases}
C_1(0) =  \frac14 \left(-\alpha_1+\frac{\alpha_2+\beta_2}{2} + \sqrt{(\alpha_2-\alpha_1)(\beta_2-\alpha_1)}\right)^2 - \left(\frac{\beta_2-\alpha_2}{4}\right)^2, \\
C_2(0) = \left(\frac{\beta_2-\alpha_2}{4}\right)^2,
\end{cases}
\\
\label{ODEbc2}
& \begin{cases}
C_1(1) = \left(\frac{\beta_1-\alpha_1}{4}\right)^2, \\
C_2(1) = \frac14 \left(\beta_2-\frac{\alpha_1+\beta_1}{2} + \sqrt{(\beta_2-\beta_1)(\beta_2-\alpha_1)}\right)^2 - \left(\frac{\beta_1-\alpha_1}{4}\right)^2,
\end{cases}
\end{align}
where
\begin{align*}
 A_1(s) = s^2 C_1(s), \quad
  A_2(s) = (1-s)^2 C_2(s), \quad  B_2(s)-B_1(s) = \sqrt{C_1(s) + C_2(s)},
\end{align*}
and~\eqref{eq:Bs}.
Moreover, $A_1,A_2,B_1,B_2$ are constant on the interval $[c_1,c_2]$.
\end{theorem}

\noindent
\textbf{Remarks.} \textbf{1)} We note that general (and conditional) part 1) of Theorem~\ref{thm:main} admits presence inside $[0, 1]$ of a subdomain,  where $A_1,A_2,B_1,B_2$ are constant. For Angelesco systems it is a generic situation which happens when ``pushing'' is not active, see \cite{GR81}.

%

\noindent
\textbf{2)} We note that known information about support of zero counting measure of MOP for Angelesco system (see \cite{GR81}) allows us to identify the subdomain where   $A_1,A_2,B_1,B_2$ are constant, i.e. interval $[c_1,c_2]$. Then it is possible, using BC \eqref{ODEbc} and \eqref{ODEbc2} to solve the ODE system \eqref{ODEsys} on $[0, c_1]$ and $[c_2, 1]$.


\subsection{Proof of Theorem~\ref{thm:main}}
Taking $d=2$ in ~\eqref{PDE1} (with $i=1, j=2$),  ~\eqref{PDE2} (with $i=1, j=2$), and  ~\eqref{PDE3} (with $i=2, j=1$; then $i=1, j=2$) gives us four ODE's:
\begin{align}
\label{ode3} & B_1'(s) s + B_2'(s)(1-s) = 0; \\
\label{ode4} & B_1(s) B_1'(s) s +B_2(s)B_2'(s) (1-s) + A_1'(s)+A_2'(s)=0;\\
\label{ode1}
& A_1(s)(B_1'(s)-B_2'(s))(1-s) + A_1'(s)(B_1(s) -B_2(s))s =0; \\
\label{ode2} & A_2(s)(B_1'(s)-B_2'(s))s + A_2'(s)(B_1(s) -B_2(s))(1-s) =0.
\end{align}
Let us simplify this system. First of all, let
$$
B(s) = B_2(s) - B_1(s).
$$
Using~\eqref{ode3} and~\eqref{ode4}, we get $B_1' = \frac{A_1'+A_2'}{sB}$, $B_2' = -\frac{A_1'+A_2'}{(1-s)B}$, so $B' = B_2'-B_1' =  -\frac{A_1'+A_2'}{s(1-s)B}$. This equation together~\eqref{ode1} and~\eqref{ode2} established~\eqref{3ODE}.
Part 1) of the Theorem~\ref{thm:main} is proved.

\medskip

Let us divide interval $[0,1]$ into two disjoint sets:
$$
I_1=\overline{\{s\in[0,1]: A_1'(s)=A_2'(s)=B'(s)=0\}} \quad \text{ and } \quad I_2=[0,1]\setminus I_1.
$$
From \cite{GR81} we know
  that: $I_1$
consists of one point if $\Delta_1$ and $\Delta_2$ are touching, and otherwise $I_1$ is an interval $[c_1,c_2]$ inside $(0,1)$.

For $s\in I_2$, the determinant of the matrix in~\eqref{3ODE} must be zero, i.e.,
\begin{equation}
s(1-s)B(s)^3 - \tfrac{1-s}{s} A_1(s) B(s) - \tfrac{s}{1-s} A_2(s) B(s) =0,
\end{equation}
which implies
\begin{equation}\label{eq:bSquared}
B(s)^2=  \tfrac{1}{s^2} A_1(s) + \tfrac{1}{(1-s)^2} A_2(s)
\end{equation}
on the set where $B(s)\ne 0$. This means that
$$
2B(s)B'(s) = \tfrac{1}{s^2} A_1'(s) + \tfrac{1}{(1-s)^2} A_2'(s) -\tfrac{2}{s^3} A_1(s) + \tfrac{2}{(1-s)^3} A_2(s).
$$
Plugging this into the third equation of ~\eqref{3ODE}, we get
$$
\tfrac{2}{s(1-s)}  (A_1'(s)+A_2'(s))+\tfrac{1}{s^2} A_1'(s) + \tfrac{1}{(1-s)^2} A_2'(s) -\tfrac{2}{s^3} A_1(s) + \tfrac{2}{(1-s)^3} A_2(s)  =0,
$$
which simplifies to
\begin{equation}
\tfrac{1+s}{s}  A_1'(s)+ \tfrac{2-s}{1-s} A_2'(s) -\tfrac{2(1-s)}{s^2} A_1(s) + \tfrac{2s}{(1-s)^2} A_2(s)  =0.
\end{equation}
The first two equations in~\eqref{3ODE} can be solved for $\tfrac{B'(s)}{B(s)}$ giving us
\begin{equation}\label{eq:As}
\tfrac{s}{1-s} \tfrac{A_1'(s)}{A_1(s)} = \tfrac{1-s}{s} \tfrac{A_2'(s)}{A_2(s)}.
\end{equation}
So our new system of two ODE's is
\begin{align}
\label{newODE1} & \tfrac{1+s}{s}  A_1'(s)+ \tfrac{2-s}{1-s} A_2'(s) -\tfrac{2(1-s)}{s^2} A_1(s) + \tfrac{2s}{(1-s)^2} A_2(s)  =0,
\\
\label{newODE2} & \tfrac{s}{1-s} \tfrac{A_1'(s)}{A_1(s)} = \tfrac{1-s}{s} \tfrac{A_2'(s)}{A_2(s)}
\end{align}
for $s\in I_2$.

%

It is not hard to notice from~\eqref{eq:bSquared} that $A_1$ and $A_2$ have double zeros at $0$ and $1$, respectively. So let
\begin{align*}
 C_1(s) = \frac{A_1(s)}{s^2}, \qquad
 C_2(s) = \frac{A_2(s)}{(1-s)^2}.
\end{align*}

Then our system~\eqref{3ODE} becomes:
\begin{equation}\label{3ODE4}
\begin{pmatrix}
-s B(s) & 0 & (1-s)C_1(s) \\
0 & -(1-s) B(s) & s C_2(s) \\
\tfrac{s}{(1-s)} & \tfrac{1-s}{s} & -B(s)
\end{pmatrix}
\begin{pmatrix}
C_1'(s) \\
C_2'(s) \\
-B'(s)
\end{pmatrix}
=
\begin{pmatrix}
2B(s)C_1(s) \\
-2B(s)C_2(s) \\
-\tfrac{2}{1-s} C_1(s) +\tfrac{2}{s} C_2(s)
\end{pmatrix}
\end{equation}

Using $B(s)^2 = C_1(s) + C_2(s)$, we can eliminate $B$:
\begin{align}
\label{2ODE1} & (1+s) s C_1'(s) + 4s  C_1(s) + (2-s)(1-s) C_2'(s)   -4(1-s) C_2(s)  =0,
\\
\label{2ODE2} &  \tfrac{s^2 C_1'(s)}{C_1(s)} +2s =  \tfrac{(1-s)^2 C_2'(s)}{C_2(s)} -2(1-s).
\end{align}

Finally, let us deal with the boundary conditions for all of our functions.

 Since our system $\{\mu_j\}_{j=1}^2$ is from a multiple Nevai class, we also have that $\mu_1$ and $\mu_2$ are in the (scalar) Nevai class $N({A_1(1)},B_1(1))$ and $N({A_2(0)},B_2(0))$, respectively. Since $\supp (\mu_j)=[\alpha_j,\beta_j]$, Weyl's theorem (see Section~\ref{s:OP}) gives us: 
%
\begin{align}
\label{newBC1} & A_1(1) = \left(\frac{\beta_1-\alpha_1}{4}\right)^2 , \quad & B_1(1) = \frac{\alpha_1+\beta_1}{2}, \\
\label{newBC4} & A_2(0) = \left(\frac{\beta_2-\alpha_2}{4}\right)^2  , \quad & B_2(0) =  \frac{\alpha_2+\beta_2}{2}.
\end{align}
The marginal conditions~\eqref{eq:marginal} give us
\begin{align*}
& A_1(0) = 0, \quad & A_2(1) = 0.
\end{align*}

We also need the other two boundary conditions
\begin{align*}
& B_1(0) = \frac12 \left(\alpha_1+\frac{\alpha_2+\beta_2}{2} - \sqrt{(\alpha_2-\alpha_1)(\beta_2-\alpha_1)}\right), \\
& B_2(1) = \frac12 \left(\beta_2+\frac{\alpha_1+\beta_1}{2} + \sqrt{(\beta_2-\beta_1)(\beta_2-\alpha_1)}\right),
\end{align*}
which can be obtained from Section 4 below.

This means that $B$ has boundary values
\begin{align*}
& B(0) = \frac12 \left(-\alpha_1+\frac{\alpha_2+\beta_2}{2} + \sqrt{(\alpha_2-\alpha_1)(\beta_2-\alpha_1)}\right), \\
 & B(1) = \frac12 \left(\beta_2-\frac{\alpha_1+\beta_1}{2} + \sqrt{(\beta_2-\beta_1)(\beta_2-\alpha_1)}\right).
\end{align*}

Finally, to get the boundary values for $C_1(t)$ and $C_2(t)$, we recall that on the region $I_2$ where $B\ne 0$, we have
\begin{equation}\label{EasAB}
B(s)^2 = \tfrac{1}{s^2} A_1(s) + \tfrac{1}{(1-s)^2} A_2(s)
\end{equation}
Taking $s\to 0$, we therefore get $A_1(0) = A_1'(0) = 0$ and
$$
B(0)^2 = \tfrac12 A_1''(0) + A_2(0),
$$
which implies $A_1''(0) = 2 B(0)^2-2A_2(0) 
$.
Similarly, $s\to 1$ gives us $A_2(1) = A_2'(1) = 0$ and
$$
B(1)^2 = A_1(1)+ \tfrac12 A_2''(1),
$$
which implies $A_2''(1) =2B(1)^2- 2A_1(1) 
$. Then $C_1(0)= \tfrac12 A_1''(0)$, $C_2(0) = A_2(0)$, $C_1(1) = A_1(1)$, $C_2(1) = \tfrac12 A_2''(1)$ which result in our boundary conditions~\eqref{ODEbc}--\eqref{ODEbc2}.
Part 2) of the theorem is now proved.

\section{Determination of the limits by means of  parametrization of $\mathfrak{R}_{\vec{t}}$}\label{AngParametriz}


In this section we employ an algebraically-geometric approach in order to
determine the limits of the NNRR's coefficients. We restrict the
consideration to the case of Angelesco system with two orthogonality measures (we allow the supports to have a common endpoint). Thus in this
setting we set
$$
d=2, \qquad \vec{t}=(t_{1},t_{2}), \quad |\vec{t}|=t_{1}+t_{2}=1, \qquad \vec{s}=s=t_{1}\in(0,1).
$$
Our input is the supports measures of
orthogonality \eqref{angelesco}
\begin{equation} \label{angelescoV}
 [\alpha_i,\beta_i] ,\,\,\,i=1,2,\quad \text{with } \quad \alpha_1<\beta_1\leqslant\alpha_{2}<\beta_2.
\end{equation}
Note that using the linear map
$y(x)=(x-\beta_{1})/(\beta_{2}-\beta_{1})$, these segments can be
transformed to
\begin{equation} \label{Star}
 [-\alpha,\,0], \quad [\beta,\,1], \qquad \alpha>0,\,\,
\beta\in [0,1),
\end{equation}
where $y(\alpha_{1})=- \alpha$ and $y(\alpha_{2})=\beta$.
Thus, without  loss of
generality, we can use \eqref{Star} as the input.

\smallskip

Our goal is to construct the following procedure: based on Theorem~\ref{thm:recurrence}, 
find the limits
\eqref{Afunc}, \eqref{Bfunc} via computing the
residues of $\Upsilon_i,\,\,i=1,2$.

\medskip

In order to reach this goal we have to solve two
problems:

\smallskip

\textbf{Problem 1.} For each $s\in(0,1)$, find the segments $[\alpha_{s,i},\beta_{s,i}]
,\,\,\,i=1,2$, of the support of the extremal vector-measure
$\vec{\omega}=(\omega_{1}, \omega_{2})$, minimizing the  energy
functional \eqref{ExtrPr}.

\smallskip

\textbf{Problem 2.} Using the endpoints
$\{\alpha_{s,i},\beta_{s,i}\}_{i=1}^{2}$ as the branch points of the
Riemann surface~$\mathfrak{R}_s$ ({which is defined in
subsection~\ref{susec1.3}}), find the limits $\{A_{s,i}, B_{s,i}\}_{i=1}^{2}$ by computing the residues of the meromorphic on
$\mathfrak{R}_s$ functions $\Upsilon_i,\,\,i=1,2$.

\smallskip

\subsection{Parametrization
of $\mathfrak{R}_{s}$ and solution to Problem 2}\label{ParRSandPr2}
To solve both problems  we use (introduced in \cite{ApKaLyTu09}
and developed in \cite{LyTu17}, \cite{LyTu}) parametrization of the
three-sheeted Riemann surfaces with four branch points.

\medskip

We fix $s\in (0, 1)$ and start with parametrization of
$\mathfrak{R}(\alpha, \beta):=\mathfrak{R}_s$, where we take \eqref{Star} for the intervals
$[\alpha_{s,i},\beta_{s,i}], i=1,2$. We
define
\begin{equation}\label{Diez}
\mathcal{U}(u):=\frac{u(2-u)^3}{(2u-1)^3},\quad u \in (1,2),\qquad \quad
{R}_{u}(\tau):=\frac{\tau^{2}(\tau+u-2)}{(2u-1)\tau-u},\quad \tau \in
\overline{\mathbb{C}}.
\end{equation}
It is not difficult to check that for $\alpha, \beta$ given in
\eqref{Star} there exists a unique solution of the equation
\begin{equation}\label{Tri1}
\exists \, \mbox{!} \, u_{\alpha, \beta}\in (1,2)\,\,\mbox{:} \qquad
\mathcal{U}(u_{\alpha, \beta})=\frac{\beta(1+\alpha)}{{\alpha + \beta}}.
\end{equation}
We have the following
\begin{theorem}[\cite{ApKaLyTu09, LyTu17, LyTu}]
\label{thm:paramRS} Riemann surface $\mathfrak{R}(\alpha,\beta)$
can be defined by means of the conformal map of the Riemann sphere $\overline{\mathbb{C}} \ni w \mapsto \textbf{z}(w) \in \mathfrak{R}(\alpha,\beta)$ given by
\begin{equation}\label{0.1}
z(w):=\pi(\mbox{\textbf{z}}(w))=
\dfrac{\alpha R_{u_{\alpha\beta}}(w)}{1+\alpha-R_{u_{\alpha\beta}}(w)}\, ,
\end{equation}
where
$\pi:\mathfrak{R}(\alpha,\beta)\to \overline{\mathbb{C}}$ is
the natural projection.
\end{theorem}

Let $\tau_0, \tau_1, \tau_2$ be $\textbf{z}^{-1}(\infty^{(0)}), \textbf{z}^{-1}(\infty^{(1)}), \textbf{z}^{-1}(\infty^{(2)})$, respectively.
Substituting $R_u$ from \eqref{Diez} into \eqref{0.1}, we obtain
\begin{equation}\label{0.10}
z(w)=-\dfrac{\alpha w^2(w-\gamma)}{(w-\tau_0)(w-\tau_1)(w-\tau_2)}\,,\qquad \gamma=2-u_{\alpha\beta},
\end{equation}
where $\tau_0\equiv \tau_{\alpha\beta}\;$
satisfies
\begin{equation}\label{0.2}
\exists\,!\, \tau_{\alpha\beta}>1:\quad 1+\alpha=R_{u_{\alpha\beta}}(\tau_{\alpha\beta})
\end{equation}
and $\tau_1, \tau_2$ are roots of the quadratic equation
\begin{equation}\label{0.3}
\tau_1+\tau_2=-(u_{\alpha\beta}+\tau_0-2),\quad \tau_1\tau_2=
-\dfrac{u_{\alpha\beta}\tau_0(u_{\alpha\beta}+\tau_0-2)}
{2u_{\alpha\beta}\tau_0-u_{\alpha\beta}-\tau_0},\quad\tau_1<\tau_2<\tau_0.
\end{equation}

\medskip

Solution of {Problem 2} is given by the
following corollary of Theorem ~\ref{thm:paramRS}.
\begin{corollary}
\label{Cor:Thm5} Let \eqref{Star} be supports \eqref{supports} of
extremal measures \eqref{ExtrPr} for some fixed~$s\in(0,1)$ of
Angelesco system \eqref{angelescoV}, and let $(u_{\alpha\beta},
\tau_{\alpha\beta})$ be the images of transformations \eqref{Tri1},
\eqref{0.2}. Then for limits \eqref{limit} of the corresponding
NNRR coefficients we have
\begin{equation}\label{0.18}
A_1(s)=-\dfrac{\alpha\tau_0^2\,C_1\,(\tau_0-\gamma)}{(\tau_0-\tau_1)^2(\tau_0-\tau_2)}, \qquad
B_1(s)=\dfrac{\alpha\tau_0\,D_{1}}{(\tau_0-\tau_1)^2
(\tau_0-\tau_2)^2},
\end{equation}
where parameters $\tau_0:=\tau_{\alpha\beta}, \tau_1, \tau_2$ are
defined in \eqref{0.2}, \eqref{0.3}, and
\begin{equation}\label{0.17}
C_1:=-\dfrac{\alpha\tau_1^2(\tau_1-\gamma)}{(\tau_0-\tau_1)^2(\tau_1-\tau_2)},\qquad
D_{1}:=\tau_0^2\tau_2+2\tau_0^2\tau_1-3\tau_0\tau_1\tau_2-\gamma\tau_0^2-\gamma\tau_1\tau_0+2\gamma\tau_1\tau_2.
\end{equation}
Formulas for  $A_2,B_2$ can be obtained by the swap of indices  $1\to2,
2\to1$.
\end{corollary}
Proof of this corollary is presented below in
subsection~\ref{ProofCol}.

\medskip

\subsection{Parametrization of supports and 
ray directions and solution to Problem~1}\label{ParOmandPr1}
Before we start dealing with {Problem 1}, let us come back to the parametrization
\eqref{Diez} and consider $(u,\tau)$ on the half-strip $
\bigsqcup:=(1,2)\times(1,\infty)$. If we invert map
\eqref{Tri1}, \eqref{0.2}, then we get a  smooth diffeomorphism
$(\mathcal{A},
\mathcal{B}):\;\bigsqcup \ni (u,\tau) \mapsto (\alpha,\beta)\in (0,+\infty) \times(0,1)$:
\begin{equation}\label{0.4}
\mathcal{A}(u,\tau):=R_u(\tau)-1,\quad
\mathcal{B}(u,\tau):=\dfrac{\mathcal{A}(u,\tau)\,\mathcal{U}(u)}
{1+\mathcal{A}(u,\tau)-\mathcal{U}(u)},
\end{equation}
which by means of coordinates $(u,\tau)$ parametrize the branch points
$\{-\alpha,\beta\}$ of the Riemann surface
$\mathfrak{R}(\alpha,\beta)$, i.e., the left endpoints of segments
\eqref{Star} of supports of the extremal vector-measure
$\overrightarrow{\omega}$.

In \cite{LyTu}  there was introduced a parametrization of the direction
$(t_1,t_2)$, see ~\eqref{multi-indices}, that corresponds to the masses for the extremal measures $(\omega_1,\omega_2)$ which have
supports $[-\alpha,0], [\beta,1]$ when the vector equilibrium problem is formulated on $[-\alpha,0], [0,1]$. It is given by the function
$$
\Theta:\;\bigsqcup\ni (u,\tau) \mapsto \theta\in (-1,1), \quad t_1=s=\frac{1+\theta}{2},\;t_2=\frac{1-\theta}{2},
$$
\begin{equation}\label{0.5}
\Theta(u,\tau):=(\tau-u)\left(\dfrac{2+2u\tau-u-\tau}
{(2u\tau-u-\tau)(u+\tau)(u+\tau-2)}\right)^{1/2}.
\end{equation}

Now we can deal with {Problem 1}. Without loss of generality (we make it
clear below in subsection~\ref{RemPr1NonT}), it is enough to
consider the Angelesco system on touching intervals ($\beta=0$):
\begin{equation}\label{0.6}
[-\alpha,0],\qquad [0,1]; \qquad  \quad \alpha>0.
\end{equation}
{Problem 1} can be decomposed into two parts:

\medskip

\textbf{Problem 1.1.} Given $\alpha$, find $s_\alpha\in(0,1)$ such
that segments \eqref{0.6} are supports of the extremal measure
of problem \eqref{ExtrPr}.

\smallskip

\textbf{Problem 1.2.}  For fixed $s\in(s_\alpha,1)$ find  the value
of $\beta_s$ so that:
\begin{equation}\label{0.7}
\supp\,\omega_1=[-\alpha,0],\quad \supp\,\omega_2=[\beta_s,1].
\end{equation}
Solution of these problems is given in the following theorem.

\begin{theorem}[for proof see \cite{ LyTu}]
\label{thm:supp-param} Given $\alpha$ in \eqref{0.6}:

 1) Excluding variable $\tau$ from the system
of equations
$$
\left\{
\begin{array}{l}
\mathcal{A}(2,{\tau})=\alpha\smallskip\\
{\Theta}(2,{\tau})=\theta
\end{array}
\right.
$$
we get the value of $\theta(\alpha)=:\theta_\alpha$. Then the
answer to {Problem 1.1} is $s_\alpha =
\displaystyle\frac{1+\theta_\alpha}{2}$.

\medskip

2) For each $s\in(s_\alpha,1)$, let $\theta = 2s-1 \in (\theta_\alpha,1)$. Then
the system
$$
\left\{
\begin{array}{l}
\mathcal{A}(\tilde{u},\tilde{\tau})=\alpha 
\smallskip\\
{\Theta}(\tilde{u},\tilde{\tau})=\theta 
\end{array}
\right.
$$
has a unique solution $(\tilde{u},\tilde{\tau})$,
and $\beta_{s}:=\mathcal{B}(\tilde{u},\tilde{\tau})$ is the
answer to {Problem 1.2.}
\end{theorem}

Summarizing, we have for $d=2$ the following \textbf{Procedure} for
finding  limits \eqref{Afunc}, \eqref{Bfunc} of NNRR coefficients
$$ A_j(s)=\lim\limits_{\mathcal N} a_{\vec{n},j},\quad
B_j(s)=\lim\limits_{\mathcal N} b_{\vec{n},j},\quad j=1,2 $$ for
the Angelesco systems of MOPs \eqref{ortho}, \eqref{angelesco}
defined on  intervals \eqref{0.6}.

\medskip

\textbf{1.} Solve {Problem 1.1:} find $\theta_\alpha, s_\alpha=
(1+\theta_\alpha)/2$. To do this, evaluate
functions $\mathcal{A}(2,\tau), \Theta(2,\tau)$ by \eqref{0.4}, \eqref{0.5}, which determines the value and $\theta_\alpha$ according to
Theorem~\ref{thm:supp-param} 1).

\smallskip

\textbf{2.} For each $s\in (s_\alpha,1)$ solve {Problem 1.2:}
find $\beta_s$ from \eqref{0.7}. To do this, solve the system from
Theorem~\ref{thm:supp-param} 2) for $\theta:=2s-1$ and substitute its
solution $(\tilde{u},\tilde{\tau})$  into the function
$\mathcal{B}$ to find $\beta_s$.

\smallskip

\textbf{3.} For each $s\in (s_\alpha,1)$ find $A_j(s), B_j(s),\;j=1,2$.
To do this, apply Corollary of Theorem~\ref{thm:paramRS} with the supports of  the
extremal measure being $[-\alpha,0]$ and $[\beta_s,1]$, i.e.,
solve equations \eqref{Tri1}, \eqref{0.2}, \eqref{0.3} 
and substitute the resulting $u_{\alpha\beta},
\tau_{\alpha\beta}, \tau_1, \tau_2$ into the formulas ~\eqref{0.18} for $A_j(s),
B_j(s),\;j=1,2$.

\medskip

To find limits \eqref{Afunc}, \eqref{Bfunc} for $s\in(0,s_\alpha)$, we do the following:

\medskip

\textbf{4.} We make reflection with respect to 0 and scaling (by $k=1/\alpha$) to get the system of intervals to the form \eqref{0.6}. As a result, the new intervals are $[-\widehat{\alpha},0]$, $[0,1]$ with $\widehat{\alpha} = 1/\alpha$.

\smallskip

\textbf{5.} We apply the above steps {1, 2, 3} of the {Procedure} to this new system of intervals to get the limits $\widehat{A}_j(s)$, $\widehat{B}_j(s)$, $j=1,2$ for $s\in(\widehat{s}_{\widehat{\alpha}},1)$ (note that $\widehat{s}_{\widehat{\alpha}} = 1-s_\alpha$).

\smallskip

\textbf{6.} Then $A_j(s) = \widehat{A}(1-s)/k^2$ and $B_j(s) = - \widehat{B}(1-s)/k$ for $j=1,2$ and $s\in(0,s_\alpha)$. Indeed, scaling by $k$ stretches all the $b_{\vec n,j}$-coefficients by $k$ and all the $a_{\vec n,j}$-coefficients by $k^2$. Reflection multiplies the $b_{\vec n,j}$-coefficients by $-1$, keeps $a_{\vec n,j}$'s intact and flips $s$ to $1-s$.


\subsection{Remark on \textbf{Problem 1} for the
measures with non-touching supports}\label{RemPr1NonT}

At first we provide an equivalent characterization of the extremal
vector-measure $\vec{\omega}=(\omega_1, \omega_2)$ of the
functional \eqref{ExtrPr}. We have (see \cite{GR81}):
\begin{equation}\label{0.8}
\left\lbrace
\begin{array}{l}
2V^{\omega_1}(z)+V^{\omega_2}(z)\quad \left\lbrace
\begin{array}{l}
\geqslant\gamma_1,\;z\in[\alpha_1,\beta_1]\\
=\gamma_1,\;z\in \supp \,\omega_1\subseteq[\alpha_1,\beta_1],\\
\end{array}
\right.\quad |\omega_1|=s=\dfrac{1+\theta}{2},\\
\\
V^{\omega_1}(z)+2V^{\omega_2}(z)\quad \left\lbrace
\begin{array}{l}
\geqslant\gamma_2,\;z\in[\alpha_2,\beta_2]\\
=\gamma_2,\;z\in \supp \, \omega_2\subseteq[\alpha_2,\beta_2],\\
\end{array}
\right. \quad |\omega_2|=1-s,
\end{array}
\right.
\end{equation}
where $V^\nu(z)=-\int\log|z-x|\, d\nu(x)$ is  log-potential of
measure $\nu$.

If we consider the vector potential
$$
\overrightarrow{W}=\left(\begin{matrix}
W_1\\
W_2
\end{matrix}\right) \;:=\mathbb{A}\,
\left(\begin{matrix}
V^{\omega_1}\\
V^{\omega_2}
\end{matrix}\right),\quad
\mathbb{A}:=\left(\begin{matrix}
1 & 2\\
2 & 1
\end{matrix}\right),
$$
where $\mathbb{A}$ is called the Angelesco matrix of interaction, then
from \eqref{0.8} we can see that components of $\overrightarrow{W}$
possess the equilibrium property
$$
\left\lbrace
\begin{array}{l}
W_1\equiv\gamma_1\quad\mbox{on}\quad \supp\,
\omega_1\subseteq[\alpha_1,\beta_1],\\
\\
W_2\equiv\gamma_2\quad\mbox{on}\quad \supp\,
\omega_2\subseteq[\alpha_2,\beta_2].
\end{array}
\right.
$$
Thus the extremal measure $\vec{\omega}$ is also called the
equilibrium measure.

\bigskip

Many properties of the equilibrium measure follow from equilibrium
relations \eqref{0.8} and from the fact that log-potential is
a convex function outside of the measure support. For example, for
$d=2$ the strict inclusion
$\supp\,\omega_i\subsetneq[\alpha_i,\beta_i]$ may happen only for one
component $i=1$ or $i=2$.
As another example, if we have for fixed $s$
in \eqref{0.8}
$$
\supp\,\omega_1=[\alpha_1, \beta_1],\quad
\supp\,\omega_2\subsetneq[\alpha_2, \beta_2],
$$
then $\supp\,\omega_2=[\alpha_2^*, \beta_2]$ with $\alpha_2<\alpha_2^*$,
and for this $s$ the extremal measure $\vec{\omega}$ is the
same as for all Angelesco systems with supports
$$
[\alpha_1, \beta_1],\quad[\tilde{\alpha}_2, \beta_2],\quad
\mbox{where}\quad \tilde{\alpha}_2\in[\beta_1, \alpha_2^*].
$$
Using this property we can reduce the solution of {Problem 1} for
the Angelesco systems with non-touching supports \eqref{Star}
to the case \eqref{0.6} considered above. Indeed, for the non-touching
case we start with case \eqref{0.6} anyway, i.e., with intervals $
[-\alpha,0]$ and $[0,1]$ and perform step 1 of the above {Procedure}: find $s_{\alpha}$.
Then we perform a new step:

\smallskip

\textbf{1.5.} Find $s^{\beta}\in( s_\alpha, 1)$ such that for the
Angelesco system supported by \eqref{0.6} we have\footnote{This can
be done by executing step 3 of the {Procedure} for $s>s_{\alpha}$
until \eqref{0.89} happens.}
\begin{equation}\label{0.89}
\supp \,\omega_{1}=[-\alpha,\,0], \quad \supp \, \omega_{2}=[\beta,\,1].
\end{equation}

\medskip
We note, that the obtained $s^{\beta}$ is equal $c_{2}$ from point
2) of Theorem~\ref{thm:main}:
$$
s^{\beta}\,=\,c_{2}.
$$
Then, performing steps 2 and 3 for $s\in(c_{2}, 1)$ we obtain
$A_{j}(s), B_{j}(s), j=1,2$.

\medskip

In an analogous way we obtain value of $s=c_{1}<c_{2}$ and $A_{j}(s),
B_{j}(s), j=1,2$ for $s\in(0, c_{1})$. At the end we recall that
for $s\in(c_{1}, c_{2})$ limits $A_{j}(s), B_{j}(s), j=1,2$ remain
to be the constants.

\subsection{Proof of Corollary of Theorem~\ref{thm:paramRS}}\label{ProofCol}
From Theorem~\ref{thm:paramRS} we know that the function $\textbf{z}:\overline{\mathbb{C}}\to \mathfrak{R}_{s}:=\mathfrak{R}(\alpha,\beta)$ is a conformal map, where
\begin{equation}\label{0.11}
z(w)=\pi(\mbox{\textbf{z}}(w))=-\dfrac{\alpha w^2(w-\gamma)}{(w-\tau_0)(w-\tau_1)(w-\tau_2)},\qquad
\tau_0:=\tau_{\alpha,\beta}, \quad \gamma:=2-u_{\alpha,\beta},
\end{equation}
see \eqref{0.10}. Meromorphic on $\mathfrak{R}_{s}$ function
$\Upsilon_{1}$ is defined by its divisor and normalization:
\begin{equation}\label{0.12}
\Upsilon_{1}(\textbf{z})\,=\,
\left\lbrace
\begin{array}{ll}
O\left(\displaystyle\frac1z\right),& \mbox{as } \textbf{z}\rightarrow \infty^{(0)},\quad (\mbox{equiv., as } w\rightarrow \tau_{0}),
\\
\phantom{O(}z\phantom{0}\,,& \mbox{as } \textbf{z}\rightarrow \infty^{(1)},\quad (\mbox{equiv., as } w\rightarrow \tau_{1}).
\end{array}
\right.
\end{equation}
Our goal is to obtain  two terms of the power series expansion of
$\Upsilon_{1}(\textbf{z})$ at the point $\infty^{(0)}$, namely to
find the coefficients $A_{1}, B_{1}$ in
\begin{equation}\label{0.13}
\Upsilon_{1}(\textbf{z})\Bigr|_{\textbf{z}\to \infty^{(0)}}\,=\,
 \displaystyle\frac{A_1}z \,\Bigr(1+\frac{B_1}z+\cdots \Bigl).
\end{equation}
In coordinates $w\in \overline{\mathbb{C}}$ we have
\begin{equation}\label{0.14}
\Upsilon_{1}(\textbf{z}(w))=C_1\dfrac{w-\tau_0}{w-\tau_1},
\end{equation}
where $C_{1}$ is determined from the normalization at the point
$\infty^{(1)}$, see \eqref{0.12}:
\begin{equation}\label{0.15}
\left( \dfrac{\Upsilon_{1}(\textbf{z}(w))}{\textbf{z}(w)}\right)\Bigr|_{w=\tau_1}=1,
\end{equation}
For the coefficients $A_{1}, B_{1}$ we have from \eqref{0.13}
\begin{equation}\label{0.16}
A_1=\Bigl({z}(w)\Upsilon_{1}(\textbf{z}(w))\Bigr)\Bigr|_{w=\tau_0},\qquad
B_1=\left[{z}(w)\left(\dfrac{{z}(w)}{A_1}
\Upsilon_{1}(\textbf{z}(w))-1\right)\right]\Bigr|_{w=\tau_0}\,.
\end{equation}
Thus substituting \eqref{0.14}, \eqref{0.11} in \eqref{0.15} we
obtain $C_1$ in \eqref{0.17}:
$$
 C_1=\left(z(w)\dfrac{w-\tau_1}{w-\tau_0}\right)\Bigr|_{w=\tau_1}=\dfrac{-\alpha\tau_1^2(\tau_1-\gamma)}
{(\tau_0-\tau_1)^2(\tau_1-\tau_2)}
$$
Analogously, plugging \eqref{0.14}, \eqref{0.11} into \eqref{0.16} for $A_1$, we obtain \eqref{0.18}:
\begin{equation}\label{0.19}
A_1=C_1\dfrac{-\alpha \tau_0^2(\tau_0-\gamma)}{(\tau_0-\tau_1)^2(\tau_0-\tau_2)}\,
=\,\frac{\alpha^2\,\tau_0^2(\tau_0-\gamma)\,\tau_1^2(\tau_1-\gamma)}
{(\tau_0-\tau_1)^4(\tau_0-\tau_2)(\tau_1-\tau_2)}.
\end{equation}
and plugging \eqref{0.14},
\eqref{0.11}, \eqref{0.19} into \eqref{0.16} for $B_1$, we get:
$$
B_1=
\frac{-\alpha w^2 (w-\gamma)\left(\displaystyle\frac{w^2(w-\gamma)}{(w-\tau_1)^2(w-\tau_2)}
\frac{(\tau_0-\tau_1)^2(\tau_0-\tau_2)}{\tau_0^2(\tau_0-\gamma)}-1\right)}
{(w-\tau_0)(w-\tau_1)(w-\tau_2)}\Bigr|_{w=\tau_0}\,.
$$
Using the notation
$P(w):=\dfrac{w^2(w-\gamma)}{(w-\tau_1)^2(w-\tau_2)}$, we continue:
$$
B_1=\frac{-\alpha w^2 (w-\gamma)}{(w-\tau_1)(w-\tau_2)}
\dfrac{(\tau_0-\tau_1)^{2}(\tau_0-\tau_2)}{\tau_0^2(\tau_0-\gamma)}\cdot
\underbrace{\dfrac{P(w)-P(\tau_0)}{w-\tau_0}\Bigr|_{w=\tau_0}}_{=P'(\tau_0)}\;.
$$
To compute $P'(w)$ we use
$$
 \dfrac{P'(w)}{P(w)}=
\frac{2}{w}+\frac{1}{w-\gamma}-\frac{2}{w-\tau_1}-\frac{1}{w-\tau_2}\,=\,
\frac{w^2\tau_2+2w^2\tau_1-3w\tau_1\tau_2-\gamma w^2-\gamma\tau_1w+2\gamma\tau_1\tau_2}
{ w (-w + \gamma) (w - \tau_1) (w - \tau_2)},
$$
This allows us to arrive to \eqref{0.18}:
$$
B_1=P(\tau_0)\frac{\,-\alpha (\tau_0-\tau_1) \,D_{1}}
{ \tau_0 (-\tau_0+\gamma) (\tau_0-\tau_1) (\tau_0-\tau_2)}
=\dfrac{-\tau_0^2(\tau_0-\gamma)}{(\tau_0-\tau_1)^{2}(\tau_0-\tau_2)}\frac{\alpha (\tau_0-\tau_1)
 \,D_{1}}{ \tau_0 (-\tau_0+\gamma) (\tau_0-\tau_1) (\tau_0-\tau_2)}
.$$
$$
=\dfrac{\tau_0}{(\tau_0-\tau_1)^{2}(\tau_0-\tau_2)}\frac{\alpha
 \,D_{1}}{  (\tau_0-\tau_2)}
.$$
Corollary of Theorem~\ref{thm:paramRS} is proved.

\section{Comparing numerics: Angelesco system $d=2$}\label{Ang_d2Num}

\subsection{Numerics: two touching intervals}
For the Angelesco systems with two intervals we now have three methods of numerically
estimating the limits $A_1(s), A_2(s), B_1(s), B_2(s)$ ($0\le s \le 1$) of the NNRR's coefficients:
\begin{itemize}
\item[(i)] by computing $a_{\vec{n},j}$ and $b_{\vec{n},j}$
recursively (through~\eqref{Nabla1A}--\eqref{Nabla3A}, see~\cite{FHVA}) for
large enough $|\vec{n}|$;
\item[(ii)] through the system of ODE's in Section~\ref{DE_Ang_d2}
(namely,~\eqref{ODEsys});
\item[(iii)] through the algebraically-geometric approach
of Section~\ref{AngParametriz}.
\end{itemize}

\begin{figure}[h!]
\centering
\minipage{0.45\textwidth}
\includegraphics[width=\linewidth]{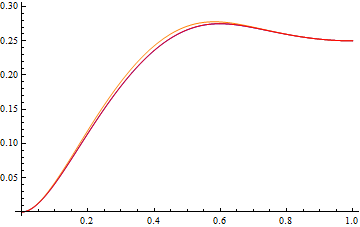}
\captionsetup{labelformat=empty}
\caption{Function $A_1(s)$}
\addtocounter{figure}{-1}
\medskip
\includegraphics[width=\linewidth]{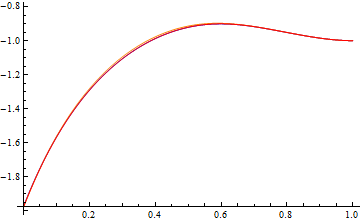}
\captionsetup{labelformat=empty}
\caption{Function $B_1(s)$}
\addtocounter{figure}{-1}
\endminipage\hfill
\minipage{0.45\textwidth}
\includegraphics[width=\linewidth]{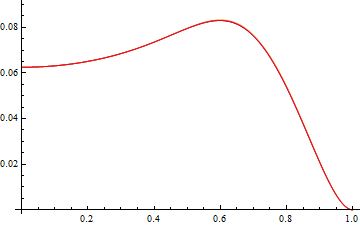}
\captionsetup{labelformat=empty}
\caption{Function $A_2(s)$}
\addtocounter{figure}{-1}
\medskip
\includegraphics[width=\linewidth]{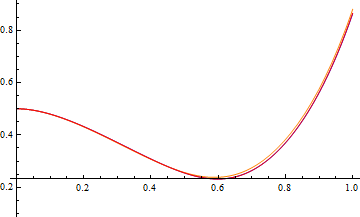}
\captionsetup{labelformat=empty}
\caption{Function $B_2(s)$}
\addtocounter{figure}{-1}
\endminipage\hfill
\caption{\small The case $\supp\,\mu_1=[-2,0]$, $\supp\,\mu_2=[0,1]$:
Blue plot: computation via the NNRR coefficients; Orange plot: computation via differential equations; Red plot: computation via the algebraically-geometric approach of Section~\ref{AngParametriz}.
}
\end{figure}
\begin{figure}[h!]
\centering
\minipage{0.45\textwidth}
\includegraphics[width=\linewidth]{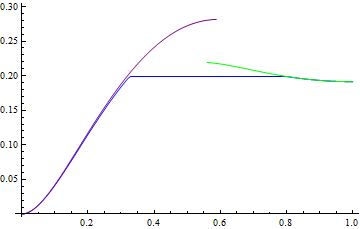}
\captionsetup{labelformat=empty}
\caption{Function $A_1(s)$}
\addtocounter{figure}{-1}
\medskip
\includegraphics[width=\linewidth]{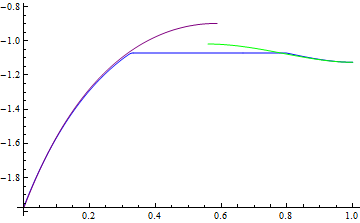}
\captionsetup{labelformat=empty}
\caption{Function $B_1(s)$}
\addtocounter{figure}{-1}
\endminipage\hfill
\minipage{0.45\textwidth}
\includegraphics[width=\linewidth]{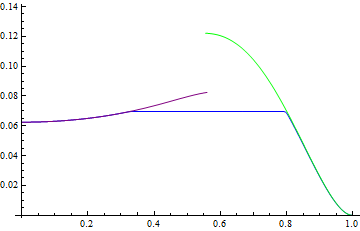}
\captionsetup{labelformat=empty}
\caption{Function $A_2(s)$}
\addtocounter{figure}{-1}
\medskip
\includegraphics[width=\linewidth]{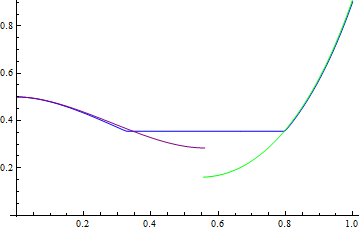}
\captionsetup{labelformat=empty}
\caption{Function $B_2(s)$}
\addtocounter{figure}{-1}
\endminipage\hfill
\caption{\small The case $\supp\,\mu_1=[-2,0]$, $\supp\,\mu_2=[0.25,1]$:
Blue plot: computation via recurrence coefficients; Purple plot: computation via differential equations with the boundary conditions at $s=0$; Green plot: computation via differential equations with the boundary conditions at $s=1$.
}
\end{figure}

On Fig.~1 we present the numerics in Wolfram Mathematica for
the case $[\alpha_1,\beta_1] = [-2,0]$, $[\alpha_2,\beta_2]=[0,1]$.
In (i) $|\vec {n}|$ was taken $1500$ (blue plot); in (ii) the in-built \texttt{NDSolve} Mathematica
 function was used (orange plot); notice that the ODE for $B_1$ in~\eqref{eq:Bs} has a singular behavior at $s=0$ and the same is true for $B_2$ at $s=1$, so one should use
 $$
 B_1'(s) = \frac{2C_1(s)+s C_1'(s)}{\sqrt{C_1(s)+C_2(s)}} \left(1+\frac{C_2(s)}{C_1(s)}\right),
 \,
 B_2'(s) = \frac{2C_2(s)-(1-s) C_2'(s)}{\sqrt{C_1(s)+C_2(s)}} \left(1+\frac{C_1(s)}{C_2(s)}\right),
 $$
 instead (these follow from ~\eqref{eq:Bs} and \eqref{eq:As}); in (iii) the interval $s\in[0,1]$ was divided into $3000$
 subintervals (red plot). The three plots are effectively indistinguishable.

\subsection{Numerics: two non-touching intervals}
On Figure 2 we present the limits $A_1(s),A_2(s),B_1(s),B_2(s)$ for an Angelesco system with $[\alpha_1,\beta_1] = [-2,0]$, $[\alpha_2,\beta_2]=[0.25,1]$. The blue plot corresponds to the computation of $a_{\vec{n},j}$ and $b_{\vec{n},j}$ recursively (via ~\eqref{Nabla1A}--\eqref{Nabla3A}) with $|\vec{n}|=1500$; the purple plot corresponds to the numerical approximation of the solution to the system of ODE's (via \eqref{ODEsys}) with the boundary conditions at $s=0$; the green plot corresponds to the numerical approximation of the solution to the system of ODE's (via \eqref{ODEsys}) with the boundary conditions at $s=1$. Equivalently, the purple plot corresponds to the coefficients' limits  for the Angelesco system with supports of $\mu_1$ and $\mu_2$ being $[-2,0.25]$ and $[0.25,1]$, while the green plot corresponds to the supports $[-2,0]$ and $[0,1]$. See Subsection~\ref{RemPr1NonT} for the explanation of this phenomenon. This can also be seen from the fact that~\eqref{ODEbc} is independent of $\beta_1$ and that ~\eqref{ODEbc2} is independent of $\alpha_2$. Again, the plots effectively overlap (away from the plateau regions).

\bigskip



\section*{Acknowledgments}

The authors are grateful to the anonymous referees for their corrections and careful proof-reading of the paper. The second author thanks W.~Van~Assche for the excellent mini-course on multiple orthogonal polynomials at the Summer School on OPSF at the University of Kent, which lead to the idea of the current paper. He also thanks the organizers of the Summer School and W.~Van~Assche and A.~Mart\'{i}nez-Finkelshtein for useful discussions.

%
%
%
%
%
%
%

\end{document}